\documentclass[reqno]{article}
\usepackage{amsmath,amsthm,amssymb}
\newtheorem{theorem}{Theorem}[section]
\newtheorem{lemma}[theorem]{Lemma}
\newtheorem{proposition}[theorem]{Proposition}
\newtheorem{corollary}[theorem]{Corollary}
\newtheorem{conjecture}[theorem]{Conjecture}
\newtheorem{problem}[theorem]{Problem}
\newtheorem{claim}{Claim}

\def\beq{\begin{equation}}\def\eeq{\end{equation}}
\def\beqn{\begin{eqnarray}}\def\eeqn{\end{eqnarray}}

\def\qed{\ifhmode\unskip\nobreak\fi\quad\ifmmode\Box\else$\Box$\fi}

\begin{document}
\title{Coverings by few monochromatic pieces - a transition between two Ramsey problems}
\author{Andr\'as Gy\'arf\'as
\\
\small Computer and Automation Research Institute\\[-0.8ex]
\small Hungarian Academy of Sciences\\[-0.8ex]
\small Budapest, P.O. Box 63\\[-0.8ex]
\small Budapest, Hungary, H-1518\\[-0.8ex]
\small \texttt{gyarfas@sztaki.hu} \and G\'{a}bor N.
S\'ark\"ozy\thanks{Research supported in part by the
National Science Foundation under Grant No. DMS-0968699.}\\
\small Computer Science Department\\[-0.8ex]
\small Worcester Polytechnic Institute\\[-0.8ex]
\small Worcester, MA, USA 01609\\[-0.8ex]
\small \texttt{gsarkozy@cs.wpi.edu}\\[-0.8ex]
\small and\\[-0.8ex]
\small Computer and Automation Research Institute\\[-0.8ex]
\small Hungarian Academy of Sciences\\[-0.8ex]
\small Budapest, P.O. Box 63\\[-0.8ex]
\small Budapest, Hungary, H-1518 \and
Stanley Selkow\\
\small Computer Science Department\\[-0.8ex]
\small Worcester Polytechnic Institute\\[-0.8ex]
\small Worcester, MA, USA 01609\\[-0.8ex]
\small \texttt{sms@cs.wpi.edu} }

\maketitle

\begin{abstract}
\footnote{Apart from a footnote referring to Pokrovskiy \cite{PO}, this manuscript was submitted to Graphs and Combinatorics in May, 2011}
The typical problem in (generalized) Ramsey theory is to find the
order of the largest monochromatic member of a family $\cal{F}$ (for
example matchings, paths, cycles, connected subgraphs) that must be
present in any edge coloring of a complete graph $K_n$ with $t$
colors. Another area is to find the minimum number of monochromatic
members of $\cal{F}$ that partition or cover the vertex set of every
edge colored complete graph. Here we propose a problem that connects
these areas: for a fixed positive integers $s\le t$, at least how
many vertices can be covered by the vertices of no more than $s$
monochromatic members of $\cal{F}$ in every edge coloring of $K_n$
with $t$ colors. Several problems and conjectures are presented,
among them a possible extension of a well-known result of Cockayne
and Lorimer on monochromatic matchings for which we prove an initial
step: every $t$-coloring of $K_n$ contains a $(t-1)$-colored
matching of size $k$ provided that $$n\ge 2k +\left\lfloor{k-1\over
2^{t-1}-1}\right\rfloor.$$

\end{abstract}

\section{Introduction}

The typical problem in (generalized) Ramsey theory is to find the
order of the largest monochromatic member of a family $\cal{F}$ (for
example matchings, paths, cycles, connected subgraphs) that must be
present in any edge coloring of a complete graph $K_n$ with $t$
colors. For easier reference these problems are called Ramsey
problems in this paper. Another well studied area, we call them {\em
cover problems}, is to find the minimum number of monochromatic
members of $\cal{F}$ that partition or cover the vertex set of every
edge colored complete graph.

Here we propose a common generalization of Ramsey and cover
problems. For a fixed positive integer $s$, at least how many
vertices can be covered by the vertices of no more than $s$
monochromatic members of $\cal{F}$ in every edge coloring of $K_n$
with $t$ colors? A somewhat related problem was proposed by Chung
and Liu \cite{CL1}: for a given graph $G$ and for fixed $s,t$, find
the smallest $n$ such that in every $t$-coloring of the edges of
$K_n$ there is a copy of $G$ colored with at most $s$ colors.

Several problems and conjectures are formulated, among them a
possible extension of a well-known result of Cockayne and Lorimer on
monochromatic matchings \cite{CL}. Our main result (Theorem \ref{main}) is
that every $t$-coloring of $K_n$ contains a $(t-1)$-colored matching
of size $k$ provided that $$n\ge 2k +\left\lfloor{k-1\over
2^{t-1}-1}\right\rfloor.$$ This result is sharp. A simple
consequence (Corollary \ref{missone}) is that every $t$-colored
$K_{2^t-2}$ has a perfect matching missing at least one color. This
is a special case of a conjecture posed in \cite{GYC}.

\subsection{Path and cycle covers}\label{path}

As far as path covers in infinite graphs are concerned, Rado
\cite{RA} has a ``perfect'' result stated here in an abridged form
with its simplified original proof.

\begin{theorem}\label{rado}
The vertex set of any $t$-colored countable complete graph can be
partitioned into finite or one-way infinite monochromatic paths,
each of a different color.
\end{theorem}

\noindent {\bf Proof.} Call a set $C\subseteq \{1,\dots,t\}$ of $k$
colors, $1\le k \le t$, {\em perfect} if there exists a set
${\cal{P}}=\{P_1,\dots,P_k\}$ of $k$ vertex disjoint finite paths
$P_1=\dots x_1,\dots,P_k=\dots,x_k$ with the following property:
$P_i$ is monochromatic in color $c_i$ and there is an infinite set
$Y$ of vertices such that $Y$ is disjoint from the paths of
$\cal{P}$ and for each $i\in\{1,\dots,k\}$ and for all $y\in Y$,
$x_iy$ is colored with $c_i$. A perfect color set exists since any
color $c_1$ present on infinitely many edges of a star incident to
vertex $x$ forms such a (one-element) set. Select a perfect set $C$
of $k$ colors so that $k$ is as large as possible ($k\le t$), this
is witnessed by $\cal{P}$ and $Y$. Let $u$ be an arbitrary vertex
not covered by $\cal{P}$. Consider a color $c$ such that $uy$ has
color $c$ for every $y\in Y^*$ where $Y^*\subseteq Y$, is infinite.
It follows from the choice of $k$ that $c\in C$. Now $u$ can be
added to the end of the $c$-colored path of $\cal{P}$, either
directly if $u\in Y$, or through a vertex $v\in Y^*$ if $u\notin Y$.
The infinite set witnessing the extension is either $Y^*$ or
$Y^*\setminus \{v\}$. Clearly the extensions can be continued to
place all vertices of the countable complete graph so that all paths
of $\cal{P}$ are finite or one-way infinite. \qed

There are several possibilities to ``finitize'' Theorem \ref{rado}.
The $2$-color version works perfectly as noted in a footnote in
\cite{GGY}.

\begin{proposition}\label{twocover}
The vertex set of any $2$-colored finite complete graph can be
partitioned into monochromatic paths, each of a different color.
\end{proposition}

\noindent {\bf Proof.} If $P_1=...,x_1,P_2=...,x_2$ are red and blue
paths and $v$ is uncovered then either $v$ can be placed as the last
vertex of one of the paths $P_i$ or one of the bypasses $P_1,x_2,v$
or $P_2,x_1,v$ extends one monochromatic path (and shortens the
other).\qed

Lehel conjectured that Proposition \ref{twocover} remains true if
paths are replaced by cycles (where the empty set, one vertex and
one edge are accepted as a cycle). Although the existence of a `near
partition' (where the two monochromatic cycles intersect in at most
one vertex) follows easily, see \cite{GY2}, it took a long time
until this was proved for large $n$ in \cite{LRSZ}, \cite{ALL}.
Recently an elementary proof was  found by Bessy and Thomass\'e
\cite{BT} that works for all $n$.

\begin{theorem}\label{twocycles}(\cite{BT})
The vertex set of any $2$-colored complete graph can be partitioned
into two monochromatic cycles of different colors.
\end{theorem}

In \cite{Gypath} and \cite{EGYP} several possible extensions of
Proposition \ref{twocover} were suggested. It turned out that for
$3$ colors one can not except full partition by distinct colors, the
first example of this phenomenon is from Heinrich \cite{HE}.
Recently the asymptotic ratios of monochromatic path and cycle
partitions with three distinct colors was obtained in \cite{GRSSz3}.

\begin{theorem}\label{3diff}(\cite{GRSSz3}) In every $3$-colored $K_n$
at least $({3\over 4}-o(1))n$ vertices can be partitioned into
monochromatic cycles of distinct colors.
\end{theorem}

We note that here the asymptotic ratio ${3\over 4}$ is best
possible. Nevertheless it was conjectured in \cite{EGYP} that
Proposition \ref{twocover} and Theorem \ref{twocycles} carries over
to any number of colors if repetitions of colors are possible.

\begin{conjecture}\label{cyccov}(\cite{EGYP}) The vertex set of every
$t$-colored complete graph can be partitioned into $t$ monochromatic
cycles.
\end{conjecture}
\footnote{Pokrovskiy \cite{PO} proved that Conjecture \label{cyccov} is not quite true for $t=3$}
The case of three colors was recently solved in asymptotic sense.

\begin{theorem}\label{3cyc}(\cite{GRSSz3}) In every
$3$-colored $K_n$ at least $(1-o(1))n$ vertices can be partitioned
into three monochromatic cycles.
\end{theorem}

The proofs of Theorems \ref{3diff}, \ref{3cyc} rely on the
Regularity Lemma and relaxations of cycles to connected matchings
(see Subsection \ref{connmatch}). The new problem we propose here is
the following.

\begin{problem}\label{spaths}
Suppose $1\le s \le t$. What is the maximum number of vertices that
can be covered by $s$ monochromatic cycles (paths) in every
$t$-coloring of the edges of $K_n$?
\end{problem}

We do not have a general conjecture here, not even for the
asymptotics (for fixed $s,t$ and large $n$). The case $s=1$ is the
path Ramsey number where the case $t=3$ (\cite{GRSSz1}) is the only
evidence that perhaps $n\over t-1$ is the true asymptotic value
(${n\over t}$ is an easy lower bound). The case $s=t$ is in
Conjecture \ref{cyccov}.

The first interesting special case is $t=3,s=2$.

\begin{conjecture}\label{conj2} In any $3$-colored $K_n$ there are two
vertex disjoint monochromatic paths (cycles) covering at least
${6n\over 7}$ vertices.
\end{conjecture}

A weaker form (for matchings instead of paths) of Conjecture
\ref{conj2} follows from Theorem \ref{main} below (when $t=3$).

\subsection{Connected matchings.}\label{connmatch}

One technique used recently in many papers (for example \cite{FL},
\cite{GRSSz1},\cite{GRSSZoneside}) in Ramsey and in covering
problems related to paths or cycles is to replace the paths or
cycles by a simpler structure, {\em monochromatic connected
matchings,} and rely on the Regularity and Blow-up Lemmas to create
paths or cycles from them. A connected monochromatic matching means
that all edges of the matching are in the same component of the
subgraph induced by the edges in the color of the matching.

Thus a natural step towards proving an asymptotic (or sharp for
large enough $n$) version of Conjecture \ref{conj2} would be to
prove it for connected matchings. However, in the problem mentioned
above we are a step behind, we could only prove Conjecture
\ref{conj2} for matchings (without the connectivity condition). A
logical plan is to treat connected pieces and matchings separately,
this is done in Subsection \ref{trees} and in Section \ref{match}.

\subsection{Covers by connected components}\label{trees}

Since every connected component contains a spanning tree of the
component, we use here the somewhat simpler tree language. A special
case of a conjecture attributed to Ryser, (appearing in his student,
Henderson's thesis \cite{HEN}) states that every intersecting
$t$-partite hypergraph has a transversal of at most $t-1$ elements.
Using the dual of the hypergraph of monochromatic components in a
$t$-coloring of complete graphs, one can easily see that the
following form of the conjecture (introduced in \cite{GY}) is
equivalent.

\begin{conjecture}\label{Rysconn} In every $t$-coloring of the edges of a complete
graph, the vertex set can be covered by the vertices of at most
$t-1$ monochromatic trees.
\end{conjecture}

If Conjecture \ref{Rysconn} is true then an easy averaging argument
would easily extend it as follows.

\begin{conjecture}\label{conncover}For every $1\le s \le t-1$ and for
every $t$-coloring of $K_n$ at least ${ns\over t-1}$ vertices can be
covered by the vertices of at most $s$ monochromatic trees.
\end{conjecture}

Since Conjecture \ref{Rysconn} is known to be true for $2\le t \le
5$, Conjecture \ref{conncover} is true for $1\le s \le t\le 5$.
Also, the case $s=1$ is known for arbitrary $t$ (originally in
\cite{GY}, \cite{GYSUR} is a recent survey). Perhaps a good test
case is to try to prove Conjecture \ref{conncover} for $s=2$ (and
for general $t$).

Since Ryser's conjecture is extended further in \cite{EGYP} by
changing cover to partition in Conjecture \ref{Rysconn}, one may
perhaps even require partition in Conjecture \ref{conncover} as
well.

\subsection{Covers by copies of a fixed graph}\label{nonbip}

It seems that to find the percentage of vertices that can be covered
by monochromatic copies of a fixed graph $H$ having at most $s$
colors is a difficult problem. Indeed, even the case when $H$ is a
single edge seems difficult. However, somewhat surprisingly, for any
fixed {\it connected non-bipartite graph} $H$ and for any fixed
$t\ge 3$ and fixed $s\le t$, the percentage of vertices of $K_n$
that can be covered by vertex disjoint monochromatic copies having
at most $s$ colors can be rather well approximated. In fact, the
following theorem can be easily obtained from the results of
\cite{LOSE}. Let $R_t(H)$ denote the smallest integer $m$ such that
in every $t$-coloring of the edges of $K_n$ there is a monochromatic
copy of $H$.

\begin{theorem}\label{nb} Suppose that $t\ge 3, 1\le s\le t$ and $H$ is a connected
non-bipartite graph. Then in every $t$-coloring of the edges of
$K_n$, at least ${s(n-R_t(H))\over t}$ vertices can be covered by
vertex disjoint monochromatic copies of $H$ using at most $s$
colors. On the other hand, for any $n$ that is divisible by $t$,
the edges of $K_n$ can be $t$-colored so that at most ${sn\over t}$
vertices can be covered by vertex disjoint monochromatic copies of
$H$ having at most $s$ colors.
\end{theorem}

\noindent \bf Proof. \rm The first part follows by selecting
successively  monochromatic copies of $H$, removing after each step
the part covered so far. Clearly, the process stops only when less
than $R_t(H)$ vertices remain. Then, an obvious averaging argument
gives that the copies in $s$ suitable colors cover the claimed
quantity.

The second part follows from the following construction. Partition
$V(K_n)$ into $t$ equal parts and color the edges within the parts
with $t$ different colors, say within part $i$ every edge gets color
$i$. The crossing edges (going from one part to another) are all
colored with the same color between any fixed pair of parts. There
are two rules. On one hand, crossing edges of color $i$ cannot be
incident to part $i$. On the other hand, the union of crossing edges
of color $i$ should span a bipartite graph. It is easy to see that
these rules can be easily guaranteed for $t\ge 3$ (and impossible to
meet for $t=2$). Because $H$ is connected, not bipartite and
crossing edges of color $i$ are not adjacent to part $i$, each
monochromatic copy of $H$ must be completely within a part.
Therefore copies of $H$ having
 at most $s$ colors are covered by the at most $s$ parts, proving the
 second statement of the theorem. \qed

\section{Covers by matchings - how to generalize Cockayne - Lorimer theorem?}\label{match}

Here we return to the most basic case, when we want to cover by
copies of an edge, i.e. by matchings. A matching in a $t$-colored
complete graph is called an $s$-colored matching if at most $s$
colors are used on its edges. To describe easily certain
$t$-colorings of $K_n$ we need, consider partition vectors with $t$
positive integer coordinates whose sum is equal to $n$. Assume that
$V(K_n)=\{1,2,\dots,n\}$. Then $[p_1,p_2,\dots,p_t]$ represents the
coloring obtained by partitioning $V(K_n)$ into parts $A_i$ so that
$|A_i|=p_i$ for $i=1,2,\dots,t$ and the color of any edge $e=(x,y)$
is the minimum $j$ for which $\{x,y\}$ has non-empty intersection
with $A_j$.

\begin{problem}\label{smatch}
Suppose $1\le s \le t$. What is the size of the largest $s$-colored
matching that can be found in every $t$-coloring of the edges of
$K_n$?
\end{problem}

The Ramsey problem, the case $s=1$ in Problem \ref{smatch}, was
completely answered by Cockayne and Lorimer \cite{CL}. Here we state
its diagonal case only.

\begin{theorem}\label{colo}(\cite{CL}) Assume $n\ge (t+1)p+2$ and
$K_n$ is arbitrarily $t$-colored. Then there is a monochromatic
matching of size $p+1$.
\end{theorem}

Observe that Theorem \ref{colo} is sharp, the coloring
$[p,p,\dots,p,2p+1]$ of $K_{(t+1)p+1}$ does not contain a
monochromatic matching with $p+1$ edges.

Notice that the case $s=t$ of Problem \ref{smatch} is trivial, any
perfect (or near-perfect if $n$ is odd) matching is obviously
optimal. In this paper we settle the case $s=t-1$, by showing that
the extremal coloring is close to the coloring
$[p,2p,4p,\dots,2^{t-1}p]$. More precisely we prove the following.

\begin{theorem}\label{main} Every $t$-coloring of $K_n$ contains a
$(t-1)$-colored matching of size $k$ provided that
$$n\ge 2k +\left\lfloor{(k-1)\over
2^{t-1}-1}\right\rfloor.$$ This is sharp for every $t\ge 2, k\ge 1$.
\end{theorem}

In case of $t=3$ Theorem \ref{main} gives Conjecture \ref{conj2} in
a weaker form.  Noting that for $k<2^{t-1}$ the second term is zero
in Theorem \ref{main}, we get the following.

\begin{corollary}\label{missone}
Every $t$-colored $K_{2^t-2}$ has a perfect matching missing at
least one color.
\end{corollary}

We note here that for $t=2,3,4$ there are results stronger than
Theorem \ref{main}. Namely, not only a $(t-1)$-colored matching of
size $k$ can be guaranteed, but a $(t-1)$-colored path on $2k$
vertices. For $t=2$ this is a well-known result \cite{GGY}, for
$t=3$ it was proved in \cite{MS} and for $t=4$ in \cite{KO}. In
fact, it was conjectured in \cite{KO} that Theorem \ref{main} holds
also if the matching of size $k$ is replaced by $P_{2k}$.

For the case  $t=4,s=2$ we suspect that the extremal coloring is
essentially $[p,p,2p,4p]$. That leads to

\begin{conjecture}\label{42}
If $n\ge \lfloor{8k-2\over 3}\rfloor$ then every $4$-coloring of
$K_n$ contains a $2$-colored matching of size $k$ .
\end{conjecture}

For the case $s=2,t=5$, the coloring $[p,p,p,2p,4p]$ and the
coloring $[p,p,p,p,2p]$ that belongs to Theorem \ref{colo} give
essentially the same parameters so we do not risk a conjecture here.
Moreover, for $s=2,t=6$ the latter coloring $[p,p,p,p,p,2p]$ is
better than $[p,p,p,p,2p,4p]$. This leads to the dilemma whether
there are better coloring in this case or $[p,p,p,p,p,2p]$ is the
extremal one? The latter possibility would be similar to the
phenomenon discussed in Subsection \ref{nonbip}, saying vaguely that
in a $6$-colored complete graph the size of the largest $2$-matching
is twice the size of the largest monochromatic matching.

\section{Large $(t-1)$-colored matchings in $t$-colored complete
graphs.}\label{proof}

Here we prove Theorem \ref{main}. To show that it is sharp, set
$N=2k-1+\left\lfloor{(k-1)\over 2^{t-1}-1}\right\rfloor=2k-1+p$,
where $p=\left\lfloor{N\over
2^t-1}\right\rfloor=\left\lfloor{k-1\over 2^{t-1}-1}\right\rfloor$.
Consider the coloring $[p,2p,4p,\dots,2^{t-2}p,q]$ of $K_N$ with
$q=N-(2^{t-1}-1)p$.

If a matching in this coloring misses color $j\ne t$ then it misses
at least $2^{j-1}p - \sum_{i<j} 2^{i-1}p=p$ vertices from the vertex
set to which color class $j$ is incident to. Thus at most $N-p=2k-1$
vertices are covered by this matching so its size is smaller than
$k$. A matching that misses color $t$ has at most $\sum_{i<t}
2^{i-1}p=(2^{t-1}-1)p\le k-1$ edges.

To prove the upper bound, consider a $t$-colored $K_n$ where
$n=2k+\left\lfloor{(k-1)\over 2^{t-1}-1}\right\rfloor=2k+p$. Set
$V=V(K_n)$, let $G_i$ denote the subgraph of $K_n$ with vertex set
$V$ and containing edges of colors different from color $i$, $1\le i
\le t$. We are going to show that for at least one $i$, $G_i$ has a
matching of size $k$. The proof is indirect: if the maximum matching
of $G_i$, $\nu(G_i)$, is at most $k-1$ for each $i$, then for the
{\em deficiency} of $G_i$, $def(G_i)$, defined as the the number of
vertices uncovered by any maximum matching of $G_i$, we have
$$def(G_i)\ge n-2\nu(G_i)\ge 2k+p-2(k-1)=p+2.$$
We apply the following well-known result, where $c_o(G)$ is the
number of odd components of $G$.

\begin{theorem}(Berge formula)\label{Berge}
$def(G)=max\{c_o(V(G)\setminus X)-|X|: X\subset V(G)\}$.
\end{theorem}

Thus, for each $i$, there is a set $X_i\subset V$ such that
\begin{equation}\label{eq}
c_o(V\setminus X_i)\ge |X_i|+p+2.
\end{equation}
Assume w.l.o.g that $|X_1|\le \dots,\le |X_t|$ and observe that the
edges between connected components of $G_i$ in $V\setminus X_i$ are
all colored with color $i$. Let $C_1^1,\dots, C_{m_1}^1$ be the
vertex sets of the connected components of $G_1$ in $V\setminus
X_1$, from (\ref{eq}) we have $m_1\ge |X_1|+p+2$.

\begin{lemma}\label{unifcov} There is  an index
$l\in \{1,2,\dots,m_1\}$, say $l=1$, such that for every $j>1$,
$\cup_{i\ne 1} C_i^1\subset X_j$.
\end{lemma}

\noindent {\bf Proof.} Suppose that $v,w\notin X_2$ where $v\in
C_q^1, w\in C_r^1$ and $q\ne r$. This implies that the edges of
color $1$ form a complete multipartite graph $M$ on $V\setminus
(X_1\cup X_2)$ with at least two partite classes. Therefore all
vertices of $M$ must be in the same connected component of $G_2$ in
$V\setminus X_2$. Thus $G_2$ has at most $1+|X_1|\le
1+|X_2|<|X_2|+p+2$ odd components in $V\setminus X_2$, contradicting
(\ref{eq}). Thus $X_2$ must cover all but at most one among the
$C_i^1$-s, say $C_1^1$ can be uncovered.

Next we show that for all $j\ge 2$, we have $\cup_{i\ne 1}
C_i^1\subset X_j$. We have seen this for $j=2$, so assume $j>2$. The
argument of the previous paragraph gives that $X_j$ covers all but
one $C_i^1$, say the exceptional one is $C_l^1$. Suppose $l\ne 1$,
say $l=2$. The inequality $|X_2|\le |X_j|$ implies $|X_2\setminus
X_j|\le |X_j\setminus X_2|$ and from this

\begin{equation}\label{eq2}
|X_1\cap (X_2\setminus X_j)|+|C_2^1\setminus X_j|\le |X_1\cap
(X_j\setminus X_2)|+|C_1^1\setminus X_2|.
\end{equation}
On the other hand, using (\ref{eq}),
$$|V\setminus X_j|=|X_1\cap (X_2\setminus X_j)|+|X_1\setminus (X_2\cup X_j)|+|C_2^1\setminus
X_j|\ge c_o(V\setminus X_j)\ge |X_j|+p+2$$
$$\ge |X_1\cap (X_j\setminus
X_2)|+|C_1^1\setminus X_2|+|X_j\cap C_2^1|+|\cup_{i\ge 3}
C_i^1|+p+2$$which can be rearranged as

\begin{equation}\label{eq3}
(|X_1\cap (X_2\setminus X_j)|+|C_2^1\setminus X_j|)-(|X_1\cap
(X_j\setminus X_2)|+|C_1^1\setminus X_2|)+|X_1\setminus (X_2\cup
X_j)|\ge $$ $$\ge |X_j\cap C_2^1|+|\cup_{i\ge 3} C_i^1|+p+2.
\end{equation}
Note that from (\ref{eq2})the left hand side of (\ref{eq3}) is at
most $|(X_2\cup X_j)|$. Thus from (\ref{eq3}) we get
$$|X_1|\ge |X_1\setminus (X_2\cup X_j)|\ge |X_j\cap
C_2^1|+|\cup_{i\ge 3} C_i^1|+p+2\ge m_1-2+p+2=m_1+p$$ and this
contradicts $m_1\ge |X_1|+p+2$ and finishes the proof of the lemma.
\qed

Call $Ker_1=C_2^1\cup \ldots \cup C_{m_1}^1$ the first kernel. With
this notation Lemma \ref{unifcov} claims that each $X_j$ with $j>1$
contains $Ker_1$. We may iterate Lemma \ref{unifcov} to define the
set $Ker_i=C_2^i\cup \ldots \cup C_{m_i}^i$, the $i$-th kernel, so
that each $X_j$ with $j>i$ contains $Ker_i$. Furthermore, these
kernels are disjoint, since $X_{i+1}$ contains $Ker_i$, but
$Ker_{i+1}$ is contained in $V\setminus X_{i+1}$. This implies that
we have the following recursion on the sizes of the $X_i$'s.

\begin{claim}\label{X_i}
For every $2\leq i \leq t$ we have
$$|X_i| \geq |X_1| + \ldots + |X_{i-1}| + (i-1)(p+1).$$
\end{claim}

Indeed, $X_i$ contains all the disjoint kernels $Ker_1, \ldots ,
Ker_{i-1}$ and thus using (\ref{eq}) we get
$$|X_i| \geq \sum_{j=1}^{i-1} |Ker_j|  \geq \sum_{j=1}^{i-1} (c_0(V\setminus X_j)-1) \geq
\sum_{j=1}^{i-1} (|X_j|+p+1),$$ as desired.

Claim \ref{X_i} implies easily by induction the following
\begin{equation}\label{|X_i|} |X_i| \geq (2^{i-1} - 1) (p+1).
\end{equation}

But then, since the kernels are disjoint, using (\ref{eq}) again we
get the following contradiction
$$n \geq \sum_{i=1}^{t} |Ker_i|  \geq \sum_{i=1}^{t} (c_0(V\setminus X_i)-1) \geq
\sum_{i=1}^{t} (|X_i|+p+1)=$$
$$=\sum_{i=1}^{t} |X_i|+t(p+1)\geq \sum_{i=1}^{t} 2^{i-1}(p+1)=
(2^t-1)(p+1) > n.$$ Here for the last inequality we have to check
$$(2^t-1)(p+1)\geq n+1=2k+p+1.$$
This is equivalent to \begin{equation}\label{p} p+1\geq
\frac{k}{2^{t-1}-1},\end{equation} which is always true for our
choice $p=\lfloor \frac{k-1}{2^{t-1}-1}\rfloor$. Indeed,
$$\left\lfloor \frac{k-1}{2^{t-1}-1}\right\rfloor = \left\lfloor
\frac{k}{2^{t-1}-1}\right\rfloor$$ (and so (\ref{p}) is trivially
true) for all cases except when $\frac{k}{2^{t-1}-1}$ is an integer,
but (\ref{p}) is true in this case as well, finishing the proof of
Theorem \ref{main}. $\Box$


\begin{thebibliography}{99}

\bibitem{ALL}P. Allen, Covering two-edge-coloured complete graphs
with two disjoint monochromatic cycles, {\em Combinatorics,
Probability and Computing}, {\bf 17} (2008), 471-486.



\bibitem{BES} S. A. Burr, P. Erd\H os, J. H. Spencer, Ramsey
theorems for multiple copies graphs, {\em Transactions of the
American Mathematical Society} {\bf 209} (1975), pp. 87-99.

\bibitem{BT} S. Bessy, S. Thomass\'e, Partitioning a graph into a
cycle and an anticycle, a proof of Lehel's conjecture, {\em Journal
of Combinatorial Theory B.}, {\bf 100} (2009), 176-180.

\bibitem{CL1} K.M. Chung, C.L. Liu, A generalization of Ramsey theory
for graphs, {\em Discrete Mathematics} {\bf 2} (1978) 117-127.

\bibitem{CL} E. J. Cockayne, P. J. Lorimer, The Ramsey number for
stripes, {\em J. Austral. Math. Soc.} {\bf 19} (1975), pp.
252-256.



\bibitem{EGYP} P. Erd\H os, A. Gy\'arf\'as, L. Pyber, Vertex
coverings by monochromatic cycles and trees, {\em Journal of
Combinatorial Theory B} {\bf 51} (1991) 90-95.

\bibitem{FL} A. Figaj, T. Luczak, The Ramsey number for a triple of
long even cycles, to appear in the {\em Journal of Combinatorial
Theory, Ser. B} {\bf 97} (2007) 584-596.



\bibitem{GGY} L. Gerencs\'er, A. Gy\'arf\'as, On Ramsey type problems, {\em Ann. Univ. Sci.
E\"otv\"os, Budapest} {\bf 10} (1967) 167 - 170.




\bibitem{GY} A. Gy\'arf\'as, Partition coverings and blocking sets
in hypergraphs (in Hungarian) {\em Communications of the Computer
and Automation Institute of the Hungarian Academy of Sciences} {\bf
71} (1977) 62 pp.


\bibitem{GY2} A. Gy\'arf\'as, Vertex coverings by monochromatic
paths and cycles, {\em Journal of Graph Theory} {\bf 7.}
(1983)131-135.

\bibitem{GYC} A. Gy\'arf\'as, Large matchings with few colors - a problem for `Eml\'ekt\'abla'
workshop, manuscript, 2010.

\bibitem{GRSSz1} A. Gy\'arf\'as, M. Ruszink\'o, G. N. S\'ark\"ozy,
E. Szemer\'edi, Three-color Ramsey numbers for paths, {\em
Combinatorica}, {\bf 27} (2007), pp. 35-69.

\bibitem{GRSSz3} A. Gy\'arf\'as, M. Ruszink\'o, G. N. S\'ark\"ozy,
E. Szemer\'edi, Partitioning 3-colored complete graphs into three
monochromatic cycles, {\em Electronic J. of Combinatorics} {\bf
18}(2011) N.53.

\bibitem{GYSSZdiam}  A. Gy\'arf\'as, G. N. S\'ark\"ozy,
E. Szemer\'edi, The Ramsey number of diamond matchings and loose
cycles in hypergraphs, {\em Electronic Journal of Combinatorics},
{\bf 15} (2008), R126.

\bibitem{GRSSZoneside} A. Gy\'arf\'as, M. Ruszink\'o, G. N. S\'ark\"ozy,
E. Szemer\'edi, One-sided coverings of complete bipartite graphs, in
Algorithms and Combinatorics 26, Topics in Discrete Mathematics,
133-144.


\bibitem{Gypath} A. Gy\'arf\'as, Monochromatic path covers, {\em
Congressus Numerantium} {\bf 109} (1995), 201-202.

\bibitem{GYSUR} A. Gy\'arf\'as, Large monochromatic components in
edge colorings of graphs - a survey, {\em `Ramsey Theory Yesterday,
Today and Tomorrow, DIMACS workshop 2008} To appear in `Progress in
Mathematics' series.

\bibitem{HE} Kathy Heinrich, personal communication, 1994.

\bibitem{HEN} J. R. Henderson, Permutation Decomposition of
(0-1)-Matrices and Decomposition Transversals, Ph.D. thesis,
Caltech, 1971.

\bibitem{KO} A. Khamseh, G. R. Omidi, A generalization of Ramsey
theory for linear forests, manuscript submitted in 2010.

\bibitem{LRSZ} T. {\L}uczak, V. R\"odl, E. Szemer\'edi, Partitioning
two-colored complete graphs into two monochromatic cycles, {\em
Combinatorics, Probability and Computing}, {\bf 7} (1998), 423-436.

\bibitem{LOSE} P.J. Lorimer, R.J. Segedin, Ramsey numbers for multiple
copies of complete graphs, {\em Journal of Graph Theory}, {\bf 2}
(1978) 89-91.

\bibitem{MS} R.Meenakshi, P.S. Sundararaghavan, Generalized Ramsey
numbers for paths in $2$-chromatic graphs, {\em Internat. J. Math.
and Math. Sci.} {\bf 9} (1986) 273-276.

\bibitem{PO} A. Pokrovskiy, Partitioning edge-coloured complete graphs into monochromatic cycles and paths, arXiv:12.05.5492v1

\bibitem{RA} R. Rado, Monochromatic paths in graphs, {\em Annals of
Discrete Mathematics} {\bf 3} (1987) 89-91.
\end{thebibliography}
\end{document}